\documentclass[12pt]{article}
\usepackage{amssymb,amsmath}
\usepackage{cases}
\usepackage{amsfonts}
\usepackage{color,xcolor}
\usepackage[left=2.0cm,right=2.0cm,top=2.0cm,bottom=2.0cm]{geometry}
\usepackage[colorlinks,citecolor=blue,urlcolor=blue]{hyperref}
\usepackage[pagewise]{lineno}\modulolinenumbers

\newtheorem{theorem}{Theorem}[section]

\newtheorem{lemma}{Lemma}[section]
\newtheorem{proposition}{Proposition}[section]

\newtheorem{definition}{Definition}[section]
\newtheorem{remark}{Remark}[section]

\begin{document}

\title{Non-uniform dependence on initial data for the generalized Camassa-Holm-Novikov equation in Besov space}

\author{Xing Wu$^{1, }$\thanks{Corresponding author. ny2008wx@163.com (Xing Wu)}\;, Yanghai Yu$^2$, Yu Xiao$^1$\\
\small \it$^1$College of Information and Management Science,
Henan Agricultural University,\\
\small  Zhengzhou, Henan 450002, China\\
\small \it $^2$School of Mathematics and Statistics, Anhui Normal University, Wuhu, Anhui 241002, China}

\date{}

\maketitle\noindent{\hrulefill}

{\bf Abstract:} Considered in this paper is the generalized Camassa-Holm-Novikov equation with high order nonlinearity, which unifies the Camassa-Holm and Novikov equations as special cases.
We show that the solution map of generalized Camassa-Holm-Novikov equation is not uniformly continuous on the initial data in Besov spaces $B_{p, r}^s(\mathbb{R})$ with $s>\max\{1+\frac{1}{p}, \frac{3}{2}\}$, $1\leq p, r< \infty$ as well as in critical space $B_{2, 1}^{\frac{3}{2}}(\mathbb{R}).$  Our
result covers and improves the previous work given by Li et al. \cite{Li 2020, 1Li 2020, Li 2021}(J. Differ. Equ. 269 (2020) 8686-8700; J. Math. Fluid Mech. 22 (2020) 4:50; J. Math. Fluid Mech., (2021) 23:36).


\vskip0mm\noindent{\hrulefill}

\section{Introduction}\label{sec1}
In this paper, we are concerned with the Cauchy problem for the generalized Camassa-Holm-Novikov (gCHN) equation with high order nonlinearity proposed by Anco, Silva and Freire \cite{Anco 2015} as follows:
\begin{eqnarray}\label{eq1}
        \left\{\begin{array}{ll}
          m_t+u^km_x+(k+1)u^{k-1}u_xm=0,\\
          m=u-u_{xx},\\
         u(0, x)=u_0. \end{array}\right.
        \end{eqnarray}
here $t>0,$ $x\in \mathbb{R}$ and $k\geq 1,$ $k\in \mathbb{N}.$  As shown in \cite{Anco 2015}, the gCHN equation (\ref{eq1}) admits a local conservation law, possesses multi-peakon solutions and exhibits wave breaking phenomena. The precise blow-up scenario and global strong solutions in the framework of Besov spaces have been investigated by Yan \cite{Yan 2019}.

When $k=1$, (\ref{eq1}) reduces to the classical Camassa-Holm (CH) equation
\begin{eqnarray}\label{eq2}
          u_t-u_{xxt}+3uu_x=2u_xu_{xx}+uu_{xxx},
                  \end{eqnarray}
which was originally derived as a bi-Hamiltonian system by Fokas and Fuchssteiner \cite{Fokas 1981} in the context of the KdV model and gained prominence after Camassa-Holm \cite{Camassa 1993} independently re-derived
it as an approximation to the Euler equations of hydrodynamics. Constantin and
Lannes \cite{Constantin 2009} later educed the CH equation from the water waves equations. The CH equation is completely integrable in the sense of having a Lax pair, a bi-Hamiltonian structure as well as possessing an infinity of conservation laws, and it also admits exact peakon solutions of the form $ce^{-|x-ct|}$ \cite{Fokas 1981, Camassa 1993, Constantin 1997, Constantin 2001, Constantin 2007, Constantin 2011}.

When $k=2$, (\ref{eq1}) becomes the famous Novikov equation \cite{Novikov 2009}
\begin{eqnarray}\label{eq3}
         u_t-u_{xxt}+4u^2u_x=3uu_xu_{xx}+u^2u_{xxx}.
        \end{eqnarray}
It is shown in \cite{Home 2008} that the Novikov equation  with cubic nonlinearity shares similar properties with the CH equation, such as a Lax pair in matrix form, a bi-Hamiltonian structure, infinitely many conserved quantities and peakon solutions.

Due to their abundant physical and mathematical properties, a series of achievements have been made in the study of CH and Novikov equations. The Cauchy problem of both equations in Sobolev spaces and Besov spaces have been investigated in \cite{Constantin 1998, 1Constantin 1998, Danchin 2001, Danchin 2003, Guo 2019, Li 2016} and \cite{Ni 2011, Himonas 2012, Wu 2012, Yan 2012, Wu 2013}  respectively. Although the solution maps of the CH and Novikov equations are continuous from $B_{p, r}^s(\mathbb{R})$ into $\mathcal{C}([0, T];  B_{p, r}^s(\mathbb{R}))$ for
$s>\max\{1+\frac{1}{p}, \frac{3}{2}\}$, $1\leq p\leq \infty$, $1\leq r<\infty$ or $(s, p, r)=(\frac{3}{2}, 2, 1)$, the solution maps are not uniformly continuous \cite{Li 2020, 1Li 2020, Li 2021}.

For general $k$, the well-posedness of system (\ref{eq1}) was established in Besov space $B_{p, r}^s$ with  $s>\max\{1+\frac{1}{p}, \frac{3}{2}\}$, $1\leq p, r\leq \infty$ (however, for $r=\infty$, the continuity of the data-to-solution map is established in a weaker topology)\cite{Zhao 2014}. The critical case for $(s, p, r)=(\frac{3}{2}, 2, 1)$ was solved in \cite{Chen 2015} later. However, to the best of our knowledge, whether  the solution map of the Cauchy problem (\ref{eq1}) for general $k$ depends not uniformly continuous on the initial data has not been studied yet. In this paper, we aim at showing that the solution map  of the initial problem (\ref{eq1}) is not uniformly continuous on the initial data in Besov spaces $B_{p, r}^s(\mathbb{R})$ with $s>\max\{1+\frac{1}{p}, \frac{3}{2}\}$, $1\leq p, r< \infty$ as well as in critical space $B_{2, 1}^{\frac{3}{2}}(\mathbb{R}).$

By using the Green function $G(x)=\frac{1}{2}e^{-|x|}$ and the identity $(1-\partial_x^2)^{-1}f=G\ast f$ for all $f \in L^2(\mathbb{R})$, we can express (\ref{eq1}) in the following equivalent form
\begin{eqnarray}\label{eq4}
        \left\{\begin{array}{ll}
         u_t+u^ku_x=-G\ast(\frac{2k-1}{2}u^{k-1}u_x^2+u^{k+1})_x
         -G\ast(\frac{k-1}{2}u^{k-2}u_x^3),\\
         u(0, x)=u_0. \end{array}\right.
        \end{eqnarray}
Let $\mathbf{\mathcal{P}}(D)=-\partial_x(1-\partial_x^2)^{-1}$ and $\mathbf{\mathcal{J}}(D)=-(1-\partial_x^2)^{-1},$ we can continue to rewrite (\ref{eq4}) as follows
\begin{eqnarray}\label{eq5}
        \left\{\begin{array}{ll}
         u_t+u^ku_x=\mathbf{\mathcal{P}}(D)(\frac{2k-1}{2}u^{k-1}u_x^2+u^{k+1})
         +\mathbf{\mathcal{J}}(D)(\frac{k-1}{2}u^{k-2}u_x^3),\\
         u(0, x)=u_0. \end{array}\right.
        \end{eqnarray}

 Our main result is stated as follows.
\begin{theorem}\label{the1.1} Let
\begin{equation*}
 s>\max\{1+\frac{1}{p},\; \frac{3}{2}\},\; (p, r)\in [1, \infty]\times [1, \infty)\;\;\mathrm{or}\;\;(s,p,r)=(\frac{3}{2},2,1)
\end{equation*}
The solution map $u_0\rightarrow S_t(u_0)$ of the initial value problem (\ref{eq5}) is not uniformly continuous from any bounded subset of  $B_{p, r}^s(\mathbb{R})$ into $\mathcal{C}([0, T];  B_{p, r}^s(\mathbb{R}))$. More precisely, there exist two sequences $g_n$ and $f_n$ such that
\begin{eqnarray*}
        \|f_n\|_{B_{p, r}^s}\lesssim 1,\qquad  \qquad \qquad \lim_{n\rightarrow \infty} \|g_n\|_{B_{p, r}^s}=0,
        \end{eqnarray*}
but
\begin{eqnarray*}
      \liminf_{n\rightarrow \infty} \|S_t(f_n+g_n)-S_t(f_n)\|_{B_{p, r}^s}\gtrsim t, \qquad t\in [0, T_0],
        \end{eqnarray*}
with small positive time $T_0$ for $T_0\leq T$.
\end{theorem}
\begin{remark}\label{rem1.1}
The system (\ref{eq5}) includes the Camassa-Holm equation for $k=1$ and the Novikov equation for $k=2$, respectively, therefore our result covers and improves the work given by Li et al. \cite{Li 2020, 1Li 2020, Li 2021}.
\end{remark}

\begin{remark}\label{rem1.2}
For general $k$, system (\ref{eq5}) involves $k+1$ order nonlinearity, making its qualitative analysis more complicated.
\end{remark}

{\bf Notations}:  Given a Banach space $X$, we denote the norm of a function on $X$ by $\|\|_{X}$, and \begin{eqnarray*}
\|\cdot\|_{L_T^\infty(X)}=\sup_{0\leq t\leq T}\|\cdot\|_{X}.
\end{eqnarray*}
For $\mathbf{f}=(f_1, f_2,...,f_n)\in X$,
\begin{eqnarray*}
\|\mathbf{f}\|_{X}^2=\|f_1\|_{X}^2+\|f_2\|_{X}^2+...+\|f_n\|_{X}^2.
\end{eqnarray*}
The symbol
$A\lesssim B$ means that there is a uniform positive constant $C$ independent of $A$ and $B$ such that $A\leq CB$.

\section{Littlewood-Paley analysis}\label{sec2}
\setcounter{equation}{0}
In this section, we will review the definition of Littlewood-Paley decomposition and nonhomogeneous Besov space, and then list some useful properties which will be frequently used in the sequel. For more details, the readers can refer to \cite{Bahouri 2011}.

There exists a couple of smooth functions $(\chi,\varphi)$ valued in $[0,1]$, such that $\chi$ is supported in the ball $\mathcal{B}\triangleq \{\xi\in\mathbb{R}:|\xi|\leq \frac 4 3\}$, $\varphi$ is supported in the ring $\mathcal{C}\triangleq \{\xi\in\mathbb{R}:\frac 3 4\leq|\xi|\leq \frac 8 3\}$. Moreover,
$$\forall\,\, \xi\in\mathbb{R},\,\, \chi(\xi)+{\sum\limits_{j\geq0}\varphi(2^{-j}\xi)}=1,$$
$$\forall\,\, \xi\in\mathbb{R}\setminus\{0\},\,\, {\sum\limits_{j\in \mathbb{Z}}\varphi(2^{-j}\xi)}=1,$$
$$|j-j'|\geq 2\Rightarrow\textrm{Supp}\,\ \varphi(2^{-j}\cdot)\cap \textrm{Supp}\,\, \varphi(2^{-j'}\cdot)=\emptyset,$$
$$j\geq 1\Rightarrow\textrm{Supp}\,\, \chi(\cdot)\cap \textrm{Supp}\,\, \varphi(2^{-j}\cdot)=\emptyset.$$
Then, we can define the nonhomogeneous dyadic blocks $\Delta_j$ as follows:
$$\Delta_j{u}= 0,\,\, \text{if}\,\, j\leq -2,\quad
\Delta_{-1}{u}= \chi(D)u=\mathcal{F}^{-1}(\chi \mathcal{F}u),$$
$$\Delta_j{u}= \varphi(2^{-j}D)u=\mathcal{F}^{-1}(\varphi(2^{-j}\cdot)\mathcal{F}u),\,\, \text{if} \,\, j\geq 0.$$

\begin{definition}[\cite{Bahouri 2011}]\label{de2.1}
Let $s\in\mathbb{R}$ and $1\leq p,r\leq\infty$. The nonhomogeneous Besov space $B^s_{p,r}(\mathbb{R})$ consists of all tempered distribution $u$ such that
\begin{align*}
\|u\|_{B^s_{p,r}(\mathbb{R})}\triangleq \Big|\Big|(2^{js}\|\Delta_j{u}\|_{L^p(\mathbb{R})})_{j\in \mathbb{Z}}\Big|\Big|_{\ell^r(\mathbb{Z})}<\infty.
\end{align*}
\end{definition}

In the following, we list some basic lemmas and properties about Besov space which will be frequently used in proving our main result.

\begin{lemma}[\cite{Bahouri 2011}]\label{lem2.1}
 (1) Algebraic properties: $\forall s>0,$ $B_{p, r}^s(\mathbb{R})$ $\cap$ $L^\infty(\mathbb{R})$ is a Banach algebra. $B_{p, r}^s(\mathbb{R})$ is a Banach algebra $\Leftrightarrow B_{p, r}^s(\mathbb{R})\hookrightarrow L^\infty(\mathbb{R})\Leftrightarrow s>\frac{1}{p}$ or $s=\frac{1}{p},$ $r=1$.\\
 (2) For any $s>0$ and $1\leq p,r\leq\infty$, there exists a positive constant $C=C(s,p,r)$ such that
$$\|uv\|_{B^s_{p,r}(\mathbb{R})}\leq C\Big(\|u\|_{L^{\infty}(\mathbb{R})}\|v\|_{B^s_{p,r}(\mathbb{R})}+\|v\|_{L^{\infty}(\mathbb{R})}\|u\|_{B^s_{p,r}(\mathbb{R})}\Big).$$
(3) For any $s \in \mathbb{R}$, $\mathbf{\mathcal{P}}(D)$ is continuous from $B_{p, r}^{s}$ into $B_{p, r}^{s+1}$, $\mathbf{\mathcal{J}}(D)$ is continuous from $B_{p, r}^{s}$ into $B_{p, r}^{s+2}.$\\
(4) Let  $1\leq p, r\leq \infty$ and $s>\max\{1+\frac{1}{p}, \frac{3}{2}\}$ or $s=1+\frac{1}{p}, r=1$. Then  we have
$$\|uv\|_{B_{p, r}^{s-2}(\mathbb{R})}\leq C\|u\|_{B_{p, r}^{s-2}(\mathbb{R})}\|v\|_{B_{p, r}^{s-1}(\mathbb{R})}.$$
\end{lemma}

\begin{lemma}[\cite{Bahouri 2011, Li 2017}]\label{lem2.2}
Let $1\leq p,r\leq \infty$. Assume that
\begin{eqnarray*}
\sigma> -\min\{\frac{1}{p}, 1-\frac{1}{p}\} \quad \mathrm{or}\quad \sigma> -1-\min\{\frac{1}{p}, 1-\frac{1}{p}\}\quad \mathrm{if} \quad \mathrm{div\,} v=0.
\end{eqnarray*}
There exists a constant $C=C(p,r,\sigma)$ such that for any solution to the
following linear transport equation:
\begin{equation*}
\partial_t f+v\cdot\nabla f=g,\qquad
f|_{t=0} =f_0,
\end{equation*}
the following statements hold:
\begin{align*}
\sup_{s\in [0,t]}\|f(s)\|_{B^{\sigma}_{p,r}}\leq Ce^{CV_{p}(v,t)}\Big(\|f_0\|_{B^\sigma_{p,r}}
+\int^t_0\|g(\tau)\|_{B^{\sigma}_{p,r}}d \tau\Big),
\end{align*}
with
\begin{align*}
V_{p}(v,t)=
\begin{cases}
\int_0^t \|\nabla v(s)\|_{B^{\sigma-1}_{p,r}}ds, &\quad \mathrm{if} \;\sigma>1+\frac{1}{p}\ \mathrm{or}\ \{\sigma=1+\frac{1}{p} \mbox{ and } r=1\},\\
\int_0^t \|\nabla v(s)\|_{B^{\sigma}_{p,r}}ds,&\quad\mathrm{if} \; \sigma=1+\frac{1}{p} \quad \mathrm{and} \quad r>1,\\
\int_0^t \|\nabla v(s)\|_{B^{\frac{1}{p}}_{p,\infty}\cap L^\infty}ds,&\quad\mathrm{if} \; \sigma<1+\frac{1}{p}.
\end{cases}
\end{align*}
\end{lemma}

\section{Non-uniform continuous dependence}\label{sec3}
\setcounter{equation}{0}
In this section, we will give the proof of Theorem \ref{the1.1}. Before proceeding further, we need to establish several crucial estimates to show that the solution $S_t(u_0)$ can be approximated by $u_0-t(u_0)^k\partial_xu_0+t(\mathbf{\mathcal{P}}(D)(u_0)+\mathbf{\mathcal{J}}(D)(u_0))$ in a small time  near $t=0$.

Firstly, we establish an acute estimate to bound the high order nonlinearity $\|u^k\|_{B_{p, r}^s}$ instead of $\|u\|_{B_{p, r}^s}^k.$ That is
\begin{lemma}\label{lem3.1}
 Let $k, m\in \mathbb{Z}^+$, $m\leq k,$ $1\leq p, r\leq \infty$ and $s>1+\frac{1}{p}$ or $s=1+\frac{1}{p}, r=1$. Then  we have
$$\|u^k\|_{B_{p, r}^s}\leq C\|u\|_{B_{p, r}^{s-1}}^{k-1}\|u\|_{B_{p, r}^s}\quad \mathrm{or}\quad \|u^k\|_{B_{p, r}^s}\leq C\|u\|_{L^\infty}^{k-1}\|u\|_{B_{p, r}^s},$$
and
\begin{equation*}
  \|u^{k-m}v^m\|_{B_{p, r}^s}\leq C\|u, v\|^{k-1}_{L^\infty}\|u, v\|_{B_{p, r}^s}.
\end{equation*}
\end{lemma}
{\bf Proof}\quad For $s>1+\frac{1}{p}$ or $s=1+\frac{1}{p}, r=1$, firstly using the product  law (2) and then the Banach algebra property (1) of Lemma \ref{lem2.1}, by recurrence method one has
\begin{eqnarray*}
\begin{split}
   \|u^k\|_{B_{p, r}^s}&\lesssim \|u\|_{L^\infty} \|u^{k-1}\|_{B_{p, r}^s}+\|u\|_{B_{p, r}^s}\|u^{k-1}\|_{L^\infty}\\
    &\lesssim \|u\|_{B_{p, r}^{s-1}} \|u^{k-1}\|_{B_{p, r}^s}+\|u\|_{B_{p, r}^s}\|u\|_{B_{p, r}^{s-1}}^{k-1}\\
      &\lesssim \|u\|_{B_{p, r}^{s-1}}(\|u\|_{B_{p, r}^{s-1}} \|u^{k-2}\|_{B_{p, r}^s}+\|u\|_{B_{p, r}^s}\|u\|_{B_{p, r}^{s-1}}^{k-2})+\|u\|_{B_{p, r}^s}\|u\|_{B_{p, r}^{s-1}}^{k-1}\\
       &\lesssim \|u\|_{B_{p, r}^{s-1}}^2 \|u^{k-2}\|_{B_{p, r}^s}+\|u\|_{B_{p, r}^s}\|u\|_{B_{p, r}^{s-1}}^{k-1}\\
        &\;\vdots\\
        &\lesssim \|u\|_{B_{p, r}^{s-1}}^{k-1}\|u\|_{B_{p, r}^s}.
        \end{split}
\end{eqnarray*}
The other two terms can be processed in a similar more relaxed way and they really hold for $s>0.$
Thus we finish the proof of Lemma \ref{lem3.1}.

According to Lemma \ref{lem3.1}, we have the estimates for $\mathbf{\mathcal{P}}(u)$, $\mathbf{\mathcal{J}}(u)$ and $u^ku_x$.
\begin{lemma}\label{lem3.2}
 Let $k\in \mathbb{Z}^+$, $1\leq p, r\leq \infty$ and $s>1+\frac{1}{p}$ or $s=1+\frac{1}{p}, r=1$. Then  we have
\begin{eqnarray*}
 \begin{split}
            \|\mathbf{\mathcal{P}}(D)(u)-\mathbf{\mathcal{P}}(D)(v)\|_{B_{p, r}^{s-1}}&\lesssim \|u-v\|_{B_{p, r}^{s-1}}\|u, v\|^{k-1}_{B_{p, r}^{s-1}}\|u, v\|_{B_{p, r}^s},\\
            \|\mathbf{\mathcal{J}}(D)(u)-\mathbf{\mathcal{J}}(D)(v)\|_{B_{p, r}^{s-1}}&\lesssim \|u-v\|_{B_{p, r}^{s-1}}\|u, v\|^{k-2}_{B_{p, r}^{s-1}}\|u, v\|^2_{B_{p, r}^s},\\
            \|\mathbf{\mathcal{P}}(D)(u)-\mathbf{\mathcal{P}}(D)(v)\|_{B_{p, r}^s}&\lesssim \|u-v\|_{B_{p, r}^s}\|u, v\|^{k-1}_{B_{p, r}^{s-1}}\|u, v\|_{B_{p, r}^s}+\|u-v\|_{B_{p, r}^{s-1}}\|u, v\|^{k-2}_{B_{p, r}^{s-1}}\|u, v\|^2_{B_{p, r}^s},\\
           \|\mathbf{\mathcal{J}}(D)(u)-\mathbf{\mathcal{J}}(D)(v)\|_{B_{p, r}^s}&\lesssim \|u-v\|_{B_{p, r}^{s-1}}\|u, v\|^{k-2}_{B_{p, r}^{s-1}}\|u, v\|^2_{B_{p, r}^s},\\
           \|\mathbf{\mathcal{P}}(D)(u)\|_{B_{p, r}^{s+1}}&\lesssim \|u\|^{k-1}_{B_{p, r}^{s-1}}\|u\|_{B_{p, r}^s}\|u\|_{B_{p, r}^{s+1}},\\
            \|\mathbf{\mathcal{J}}(D)(u)\|_{B_{p, r}^{s+1}}&\lesssim \|u\|^{k-1}_{B_{p, r}^{s-1}}\|u\|_{B_{p, r}^s}\|u\|_{B_{p, r}^{s+1}},\\
            \|u^ku_x\|_{B_{p, r}^{s-1}}&\lesssim \|u\|^k_{B_{p, r}^{s-1}}\|u\|_{B_{p, r}^s},\\
             \|u^ku_x\|_{B_{p, r}^{s+1}}&\lesssim \|u\|_{B_{p, r}^{s+1}}\|u\|^k_{B_{p, r}^s}+ \|u\|_{B_{p, r}^{s+2}}\|u\|^k_{B_{p, r}^{s-1}},
          \end{split}
\end{eqnarray*}
\mbox{and}
   \begin{eqnarray*}
            \|u^ku_x-v^kv_x\|_{B_{p, r}^s}\lesssim \|u-v\|_{B_{p, r}^{s-1}}\|u, v\|^{k-1}_{B_{p, r}^{s-1}}\|v\|_{B_{p, r}^{s+1}}+\|u-v\|_{B_{p, r}^s}\|u, v\|^k_{B_{p, r}^s}+\|u-v\|_{B_{p, r}^{s+1}}\|u\|^k_{B_{p, r}^{s-1}}.
   \end{eqnarray*}
\end{lemma}
{\bf Proof}\quad For $s>1+\frac{1}{p}$ or $s=1+\frac{1}{p}, r=1$, using (3), (4), (1) of Lemma \ref{lem2.1} and the Young inequality, we have
\begin{eqnarray*}
 \begin{split}
            \|\mathbf{\mathcal{P}}(D)(u)-\mathbf{\mathcal{P}}(D)(v)\|_{B_{p, r}^{s-1}}&\lesssim \|u^{k-1}u_x^2-v^{k-1}v_x^2\|_{B_{p, r}^{s-2}}+\|u^{k+1}-v^{k+1}\|_{B_{p, r}^{s-2}}\\
          &\lesssim \|u^{k-1}(u_x-v_x)(u_x+v_x)\|_{B_{p, r}^{s-2}}+\|(u^{k-1}-v^{k-1})v_x^2\|_{B_{p, r}^{s-2}}\\
          &\quad+\|u^{k+1}-v^{k+1}\|_{B_{p, r}^{s-2}}\\
          &\lesssim \|u_x-v_x\|_{B_{p, r}^{s-2}}\|u^{k-1}(u_x+v_x)\|_{B_{p, r}^{s-1}}+\|v_x^2\|_{B_{p, r}^{s-2}}\|u^{k-1}-v^{k-1}\|_{B_{p, r}^{s-1}}\\
          &\quad+\|u^{k+1}-v^{k+1}\|_{B_{p, r}^{s-2}}\\
          &\lesssim \|u-v\|_{B_{p, r}^{s-1}}\|u\|^{k-1}_{B_{p, r}^{s-1}}\|u, v\|_{B_{p, r}^s}+\|v\|_{B_{p, r}^{s-1}}\|v\|_{B_{p, r}^s}\|u-v\|_{B_{p, r}^{s-1}}\|u, v\|^{k-2}_{B_{p, r}^{s-1}}\\
           &\quad+ \|u-v\|_{B_{p, r}^{s-1}}\|u, v\|^k_{B_{p, r}^{s-1}}\\
          &\lesssim \|u-v\|_{B_{p, r}^{s-1}}\|u, v\|^{k-1}_{B_{p, r}^{s-1}}\|u, v\|_{B_{p, r}^s},\\
          \|\mathbf{\mathcal{P}}(D)(u)-\mathbf{\mathcal{P}}(D)(v)\|_{B_{p, r}^s}&\lesssim \|u^{k-1}u_x^2-v^{k-1}v_x^2\|_{B_{p, r}^{s-1}}+\|u^{k+1}-v^{k+1}\|_{B_{p, r}^{s-1}}\\
          &\lesssim \|u^{k-1}(u_x-v_x)(u_x+v_x)\|_{B_{p, r}^{s-1}}+\|(u^{k-1}-v^{k-1})v_x^2\|_{B_{p, r}^{s-1}}\\
          &\quad+\|u^{k+1}-v^{k+1}\|_{B_{p, r}^{s-1}}\\
          &\lesssim \|u\|^{k-1}_{B_{p, r}^{s-1}}\|u-v\|_{B_{p, r}^s}\|u, v\|_{B_{p, r}^s}+\|u-v\|_{B_{p, r}^{s-1}}\|u, v\|^{k-2}_{B_{p, r}^{s-1}}\|v\|^2_{B_{p, r}^s}\\
           &\quad+ \|u-v\|_{B_{p, r}^{s-1}}\|u, v\|^k_{B_{p, r}^{s-1}}\\
          &\lesssim \|u-v\|_{B_{p, r}^s}\|u, v\|^{k-1}_{B_{p, r}^{s-1}}\|u, v\|_{B_{p, r}^s}+\|u-v\|_{B_{p, r}^{s-1}}\|u, v\|^{k-2}_{B_{p, r}^{s-1}}\|u, v\|^2_{B_{p, r}^s},
          \end{split}
\end{eqnarray*}
and with the aid of the interpolation inequality, we obtain that
\begin{eqnarray*}
 \begin{split}
            \|\mathbf{\mathcal{P}}(D)(u)\|_{B_{p, r}^{s+1}}&\lesssim \|u^{k-1}u_x^2\|_{B_{p, r}^s}+\|u^{k+1}\|_{B_{p, r}^s}\\
          &\lesssim \|u^{k-1}\|_{B_{p, r}^s}\|u_x^2\|_{L^\infty}+\|u^{k-1}\|_{L^\infty} \|u_x^2\|_{B_{p, r}^s}+\|u\|^k_{B_{p, r}^{s-1}}\|u\|_{B_{p, r}^s}\\
          &\lesssim \|u\|^{k-2}_{B_{p, r}^{s-1}}\|u\|_{B_{p, r}^s}\|u\|^2_{B_{p, r}^s}+\|u\|^{k-1}_{B_{p, r}^{s-1}}\|u\|_{B_{p, r}^s}\|u\|_{B_{p, r}^{s+1}}+\|u\|^k_{B_{p, r}^{s-1}}\|u\|_{B_{p, r}^s}\\
          &\lesssim \|u\|^{k-1}_{B_{p, r}^{s-1}}\|u\|_{B_{p, r}^s}\|u\|_{B_{p, r}^{s+1}}.
          \end{split}
\end{eqnarray*}
The other terms can be processed in a similar way, here we omit the details..

We can see  the necessity of being bounded in $B_{p, r}^{s-1}$ from Lemma \ref{lem3.2}, however, for critical index $(s, p, r)=(\frac{3}{2}, 2, 1)$, there is no estimates of solutions in $B_{p, r}^{s-1}$, we can use $C^{0, 1}$ instead of $B_{p, r}^{s-1}$.

\begin{lemma}\label{lem3.3}
 Let $k\in \mathbb{Z}^+$, $(s, p, r)=(\frac{3}{2}, 2, 1)$. Then  we have
\begin{eqnarray*}
 \begin{split}
            \|\mathbf{\mathcal{P}}(D)(u)-\mathbf{\mathcal{P}}(D)(v)\|_{B_{p, r}^s}&\lesssim \|u-v\|_{B_{p, r}^s}\|u, v\|^{k-1}_{C^{0, 1}}\|u, v\|_{B_{p, r}^s},\\
           \|\mathbf{\mathcal{J}}(D)(u)-\mathbf{\mathcal{J}}(D)(v)\|_{B_{p, r}^s}&\lesssim \|u-v\|_{B_{p, r}^s}\|u, v\|^{k-1}_{C^{0, 1}}\|u, v\|_{B_{p, r}^s},\\
           \|\mathbf{\mathcal{P}}(D)(u)\|_{B_{p, r}^{s+1}}&\lesssim \|u\|^k_{C^{0, 1}}\|u\|_{B_{p, r}^{s+1}},\;\;
            \|\mathbf{\mathcal{J}}(D)(u)\|_{B_{p, r}^{s+1}}\lesssim \|u\|^k_{C^{0, 1}}\|u\|_{B_{p, r}^s},\\
            \|u^ku_x\|_{B_{p, r}^{s+1}}&\lesssim \|u\|_{B_{p, r}^{s+1}}\|u\|^k_{B_{p, r}^s}+ \|u\|_{B_{p, r}^{s+2}}\|u\|^k_{L^\infty},
          \end{split}
\end{eqnarray*}
\mbox{and}
   \begin{eqnarray*}
            \|u^ku_x-v^kv_x\|_{B_{p, r}^s}\lesssim \|u-v\|_{L^\infty}\|u, v\|^{k-1}_{C^{0, 1}}\|u, v\|_{B_{p, r}^{s+1}}+\|u-v\|_{B_{p, r}^s}\|u, v\|^k_{C^{0, 1}}+\|u-v\|_{B_{p, r}^{s+1}}\|u\|^k_{L^\infty}.
   \end{eqnarray*}
\end{lemma}

With the help of Lemmas \ref{lem3.1}-\ref{lem3.3}, we can establish the estimates of the difference between the solution $S_t(u_0)$ and initial data $u_0$ in different Besov norms.

\begin{proposition}\label{pro 3.1}
Assume that $s>\max\{1+\frac{1}{p}, \frac{3}{2}\}$, $(p, r)\in [1, \infty]\times [1, \infty)$ and $\|u_0\|_{B^s_{p,r}}\lesssim 1$. Then under the assumptions of Theorem \ref{the1.1}, we have
\begin{align*}
&\|S_{t}(u_0)-u_0\|_{B^{s-1}_{p,r}}\lesssim t\|u_0\|^{k-1}_{B^{s-1}_{p,r}}\|u_0\|^2_{B^s_{p,r}},
\\&\|S_{t}(u_0)-u_0\|_{B^{s}_{p,r}}\lesssim t(\|u_0\|^{k+1}_{B^{s}_{p,r}}+\|u_0\|^k_{B^{s-1}_{p,r}}\|u_0\|_{B^{s+1}_{p,r}}),
\\&\|S_{t}(u_0)-u_0\|_{B^{s+1}_{p,r}}\lesssim t(\|u_0\|^k_{B^{s}_{p,r}}\|u_0\|_{B^{s+1}_{p,r}}+\|u_0\|^k_{B^{s-1}_{p,r}}\|u_0\|_{B^{s+2}_{p,r}}).
\end{align*}
\end{proposition}
{\bf Proof}\quad For simplicity, denote $u(t)=S_t(u_0)$. Firstly, according to the local well-posedness result \cite{Zhao 2014}, there exists a positive time $T=T(\|u_0\|_{B^s_{p, r}}, s, p, r)$ such that the solution $u(t)$ belongs to $\mathcal{C}([0, T];  B_{p, r}^s)$. Moreover, by Lemmas \ref{lem2.1}-\ref{lem2.2}, for all $t\in[0,T]$ and $\gamma\geq s-1$,  there holds
\begin{equation}\label{eq3.1}
 \|u(t)\|_{B^\gamma_{p,r}}\leq C\|u_0\|_{B^\gamma_{p,r}}.
\end{equation}

For $t\in[0,T]$, using the differential mean value theorem, the Minkowski inequality, Lemma \ref{lem3.2} with $v=0$ and the interpolation inequality, we have from (\ref{eq3.1}) that
\begin{eqnarray*}
\begin{split}
 \|u(t)-u_0\|_{B_{p, r}^s}&\lesssim \int_0^t\|\partial_\tau u\|_{B_{p, r}^s}d\tau\\
 &\lesssim \int_0^t \|\mathbf{\mathcal{P}}(D)(u)\|_{B_{p, r}^s}d\tau+\int_0^t \|\mathbf{\mathcal{J}}(D)(u)\|_{B_{p, r}^s}d\tau+\int_0^t \|u^ku_x\|_{B_{p, r}^s}d\tau\\
 &\lesssim t(\|u\|^{k+1}_{L_t^\infty(B^{s}_{p,r})}+\|u\|^k_{L_t^\infty(B_{p, r}^{s-1})}\|u\|_{L_t^\infty(B^{s+1}_{p,r})})\\
 &\lesssim t(\|u_0\|^{k+1}_{B^{s}_{p,r}}+\|u_0\|^k_{B^{s-1}_{p,r}}\|u_0\|_{B^{s+1}_{p,r}}),\\
 \|u(t)-u_0\|_{B_{p, r}^{s-1}}&\lesssim \int_0^t\|\partial_\tau u\|_{B_{p, r}^{s-1}}d\tau\\
 &\lesssim \int_0^t \|\mathbf{\mathcal{P}}(D)(u)\|_{B_{p, r}^{s-1}}d\tau+\int_0^t \|\mathbf{\mathcal{J}}(D)(u)\|_{B_{p, r}^{s-1}}d\tau+\int_0^t \|u^ku_x\|_{B_{p, r}^{s-1}}d\tau\\
 &\lesssim t\|u\|^{k-1}_{L_t^\infty(B^{s-1}_{p,r})}\|u\|^2_{L_t^\infty(B^s_{p,r})}\\
 &\lesssim t\|u_0\|^{k-1}_{B^{s-1}_{p,r}}\|u_0\|^2_{B^s_{p,r}},
\end{split}
\end{eqnarray*}
and
\begin{eqnarray*}
\begin{split}
 \|u(t)-u_0\|_{B_{p, r}^{s+1}}&\lesssim \int_0^t\|\partial_\tau u\|_{B_{p, r}^{s+1}}d\tau\\
 &\lesssim \int_0^t \|\mathbf{\mathcal{P}}(D)(u)\|_{B_{p, r}^{s+1}}d\tau+\int_0^t \|\mathbf{\mathcal{J}}(D)(u)\|_{B_{p, r}^{s+1}}d\tau+\int_0^t \|u^ku_x\|_{B_{p, r}^{s+1}}d\tau\\
 &\lesssim t(\|u\|^k_{L_t^\infty(B^{s}_{p,r})}\|u\|_{L_t^\infty(B^{s+1}_{p,r})}
+\|u\|^k_{L_t^\infty(B^{s-1}_{p,r})}\|u\|_{L_t^\infty(B^{s+2}_{p,r})})
\\&\lesssim t(\|u_0\|^k_{B^{s}_{p,r}}\|u_0\|_{B^{s+1}_{p,r}}+\|u_0\|^k_{B^{s-1}_{p,r}}\|u_0\|_{B^{s+2}_{p,r}}).
\end{split}
\end{eqnarray*}
Thus, we finish the proof of Proposition \ref{pro 3.1}.

Since Proposition \ref{pro 3.1} fails for critical index $(s, p, r)=(\frac32, 2, 1)$ due to
the lack of the estimate of solutions in $B_{2, 1}^{\frac{1}{2}}$, it needs special treatment. That is

\begin{proposition}\label{pro 3.2}
Assume that $(s,p,r)=(\frac32,2,1)$ and $\|u_0\|_{B^s_{p,r}}\lesssim 1$. Under the assumptions of Theorem \ref{the1.1}, we have
\begin{align*}
&\|S_{t}(u_0)-u_0\|_{L^\infty}\lesssim t\|u_0\|^{k+1}_{C^{0,1}},\\
&\|S_{t}(u_0)-u_0\|_{B^{s}_{p,r}}\lesssim t(\|u_0\|^{k+1}_{B^{s}_{p,r}}+\|u_0\|^k_{C^{0, 1}}\|u_0\|_{B^{s+1}_{p,r}}),\\
&\|S_{t}(u_0)-u_0\|_{B^{s+1}_{p,r}}\lesssim t(\|u_0\|^k_{B^{s}_{p,r}}\|u_0\|_{B^{s+1}_{p,r}}+\|u_0\|^k_{C^{0, 1}}\|u_0\|_{B^{s+2}_{p,r}}),
\end{align*}
\end{proposition}
{\bf Proof}\quad For simplicity, denote $u(t)=S_t(u_0)$. According to the local well-posedness result \cite{Chen 2015}, there exists a positive time $T=T(\|u_0\|_{B^s_{p, r}}, s, p, r)$ such that the solution $u(t)$ belongs to $\mathcal{C}([0, T];  B_{p, r}^s)$. Moreover, by Lemmas \ref{lem2.1}-\ref{lem2.2} and a bootstrap from \cite{Bahouri 2011} on the Camassa-Holm equation, there holds for all $t\in[0,T]$ and $\gamma\geq \frac32$ that
\begin{align}\label{eq3.2}
\|u(t)\|_{B^\gamma_{2,1}}\leq C\|u_0\|_{B^\gamma_{2,1}}, \qquad \|u(t)\|_{C^{0,1}}\leq C\|u_0\|_{C^{0,1}}.
\end{align}
Using the equivalent form (\ref{eq4}) and convolution properties yields
\begin{eqnarray*}
\begin{split}
\|u(t)-u_0\|_{L^\infty}
&\lesssim \int^t_0\|\partial_\tau u\|_{L^\infty} d\tau \\
&\lesssim \int^t_0 \|G_x\ast(u^{k-1}u_x^2+u^{k+1})\|_{L^\infty}d\tau+\int^t_0 \|G\ast(u^{k-2}u_x^3)\|_{L^\infty}d\tau+\int^t_0 \|u^ku_x\|_{L^\infty}d\tau\\
&\lesssim t\|u\|^{k+1}_{L_t^\infty(C^{0,1})}\\
&\lesssim t\|u_0\|^{k+1}_{C^{0,1}}.
\end{split}
\end{eqnarray*}

Using Lemma \ref{lem3.3} and $B_{2, 1}^{\frac{3}{2}}(\mathbb{R})\hookrightarrow C^{0, 1}(\mathbb{R})$, following the same procedure of estimates leading to Proposition \ref{pro 3.1}, we get

\begin{eqnarray*}
\begin{split}
           \|u(t)-u_0\|_{B^{\frac{3}{2}}_{2,1}}&\lesssim \int_0^t\|\partial_\tau u\|_{B^{\frac{3}{2}}_{2,1}}d\tau\\
 &\lesssim \int_0^t \|\mathbf{\mathcal{P}}(D)(u)\|_{B^{\frac{3}{2}}_{2,1}}d\tau+\int_0^t \|\mathbf{\mathcal{J}}(D)(u)\|_{B^{\frac{3}{2}}_{2,1}}d\tau+\int_0^t \|u^ku_x\|_{B^{\frac{3}{2}}_{2,1}}d\tau\\
 &\lesssim t(\|u\|^{k+1}_{L_t^\infty(B^{\frac{3}{2}}_{2,1})}+\|u\|^k_{L_t^\infty(C^{0, 1})}
 \|u\|_{L_t^\infty(B^{\frac{5}{2}}_{2,1})})\\
           &\lesssim t(\|u_0\|^{k+1}_{B_{2, 1}^\frac{3}{2}}+\|u_0\|_{C^{0, 1}}^k\|u_0\|_{B_{2, 1}^\frac{5}{2}}),\\
           \end{split}
\end{eqnarray*}
and
\begin{eqnarray*}
\begin{split}
           \|u(t)-u_0\|_{B^{\frac{5}{2}}_{2,1}}&\lesssim \int_0^t\|\partial_\tau u\|_{B^{\frac{5}{2}}_{2,1}}d\tau\\
 &\lesssim \int_0^t \|\mathbf{\mathcal{P}}(D)(u)\|_{B^{\frac{5}{2}}_{2,1}}d\tau+\int_0^t \|\mathbf{\mathcal{J}}(D)(u)\|_{B^{\frac{5}{2}}_{2,1}}d\tau+\int_0^t \|u^ku_x\|_{B^{\frac{5}{2}}_{2,1}}d\tau\\
 &\lesssim t(\|u\|^k_{L_t^\infty(B^{\frac{3}{2}}_{2,1})}\|u\|_{L_t^\infty(B^{\frac{5}{2}}_{2,1})}
  +\|u\|^k_{L_t^\infty(C^{0, 1})}\|u\|_{L_t^\infty(B^{\frac{7}{2}}_{2,1})})\\
           &\lesssim t(\|u_0\|^k_{B_{2, 1}^\frac{3}{2}}\|u_0\|_{B_{2, 1}^\frac{5}{2}}+\|u_0\|_{C^{0, 1}}^k\|u_0\|_{B_{2, 1}^\frac{7}{2}}).
          \end{split}
\end{eqnarray*}
Thus, we complete the proof of Proposition \ref{pro 3.2}.

With the different norms  estimates  of $u-u_0$ at hand, we have the following core estimates, which implies that for specially selected initial data $u_0$ in $B^s_{p,r}$, the corresponding solution $S_t(u_0)$ can be approximated by $u_0-t(u_0)^k\partial_xu_0+t(\mathbf{\mathcal{P}}(D)(u_0)+\mathbf{\mathcal{J}}(D)(u_0))$ in a small time  near $t=0$.

\begin{proposition}\label{pro 3.3}
Assume that $\|u_0\|_{B^s_{p,r}}\lesssim 1$. Under the assumptions of Theorem \ref{the1.1}, then
\begin{itemize}
  \item for $s>\max\{1+\frac{1}{p},\frac{3}{2}\}$ and $(p, r)\in [1, \infty]\times [1, \infty),$ we have
  \begin{align*}
\|S_{t}(u_0)-u_0-t\mathbf{v}_0\|_{B^{s}_{p,r}}\lesssim t^{2}(\|u_0\|^{k+1}_{B^s_{p,r}}+\|u_0\|^k_{B^{s-1}_{p,r}}\|u_0\|_{B^{s+1}_{p,r}}
+\|u_0\|^{2k}_{B^{s-1}_{p,r}}\|u_0\|_{B^{s+2}_{p,r}}),
\end{align*}
  \item for $(s, p, r)=(\frac{3}{2}, 2, 1)$, we have
  \begin{align*}
\|S_{t}(u_0)-u_0-t\mathbf{v}_0\|_{B^{s}_{p,r}}\lesssim t^2(\|u_0\|^{k+1}_{B^s_{p,r}}+\|u_0\|^k_{C^{0, 1}}\|u_0\|_{B^{s+1}_{p,r}}
+\|u_0\|^{2k}_{C^{0, 1}}\|u_0\|_{B^{s+2}_{p,r}}).
\end{align*}
\end{itemize}
here $\mathbf{v}_0=-u_0^k\partial_x u_0+\mathbf{\mathcal{P}}(D)(u_0)+\mathbf{\mathcal{J}}(D)(u_0).$
\end{proposition}
{\bf Proof}\quad {\bf Case 1:} $s>\max\{1+\frac{1}{p},\frac{3}{2}\}$ and $(p, r)\in [1, \infty]\times [1, \infty)$.

Using Lemma \ref{lem3.2}, we obtain from  Propositions \ref{pro 3.1} and (\ref{eq3.1}) that
\begin{eqnarray*}
\begin{split}
 \|u(t)-u_0-t\mathbf{v}_0\|_{B^s_{p,r}}
&\leq \int^t_0\|\partial_\tau u-\mathbf{v}_0\|_{B^s_{p,r}} d\tau \\
&\leq \int^t_0\|\mathbf{\mathcal{P}}(D)(u)-\mathbf{\mathcal{P}}(D)(u_0)\|_{B^s_{p,r}} d\tau+\|\mathbf{\mathcal{J}}(D)(u)-\mathbf{\mathcal{J}}(D)(u_0)\|_{B^s_{p,r}} d\tau \\
&\quad+\int^t_0\|u^k\partial_xu-u^k_0\partial_xu_0\|_{B^s_{p,r}}d\tau
\\&\lesssim \int^t_0\|u(\tau)-u_0\|_{B^s_{p,r}} d\tau+\int^t_0\|u(\tau)-u_0\|_{B^{s-1}_{p,r}} \|u_0\|^{k-1}_{B^{s-1}_{p,r}} \|u_0\|_{B^{s+1}_{p,r}} d\tau\\
&\quad \ + \int^t_0\|u(\tau)-u_0\|_{B^{s+1}_{p,r}}  \|u(\tau)\|^k_{B^{s-1}_{p,r}}d\tau\\
&\lesssim t^2(\|u_0\|^{k+1}_{B^s_{p,r}}+\|u_0\|^k_{B^{s-1}_{p,r}}\|u_0\|_{B^{s+1}_{p,r}}
+\|u_0\|^{2k}_{B^{s-1}_{p,r}}\|u_0\|_{B^{s+2}_{p,r}}).
\end{split}
\end{eqnarray*}

{\bf Case 2: $(s, p, r)=(\frac{3}{2}, 2, 1)$}.

Using Lemma \ref{lem3.3}, we obtain from  Propositions \ref{pro 3.2} and (\ref{eq3.2}) that
\begin{eqnarray*}
\begin{split}
 \|u(t)-u_0-t\mathbf{v}_0\|_{B_{2, 1}^\frac{3}{2}}
&\leq \int^t_0\|\partial_\tau u-\mathbf{v}_0\|_{B_{2, 1}^\frac{3}{2}} d\tau \\
&\leq \int^t_0\|\mathbf{\mathcal{P}}(D)(u)-\mathbf{\mathcal{P}}(D)(u_0)\|_{B_{2, 1}^\frac{3}{2}} d\tau+\|\mathbf{\mathcal{J}}(D)(u)-\mathbf{\mathcal{J}}(D)(u_0)\|_{B_{2, 1}^\frac{3}{2}} d\tau \\
&\quad+\int^t_0\|u^k\partial_xu-u^k_0\partial_xu_0\|_{B_{2, 1}^\frac{3}{2}}d\tau
\\&\lesssim \int^t_0\|u(\tau)-u_0\|_{B_{2, 1}^\frac{3}{2}} d\tau+\int^t_0\|u(\tau)-u_0\|_{L^\infty} \|u_0\|^{k-1}_{C^{0, 1}} \|u_0\|_{B_{2, 1}^\frac{5}{2}} d\tau\\
&\quad \ + \int^t_0\|u(\tau)-u_0\|_{B_{2, 1}^\frac{5}{2}}  \|u(\tau)\|^k_{L^\infty}d\tau\\
&\lesssim t^2(\|u_0\|^{k+1}_{B^s_{p,r}}+\|u_0\|^k_{C^{0, 1}}\|u_0\|_{B^{s+1}_{p,r}}
+\|u_0\|^{2k}_{C^{0, 1}}\|u_0\|_{B^{s+2}_{p,r}}).
\end{split}
\end{eqnarray*}
Thus, we complete the proof of Proposition \ref{pro 3.3}.

Now, we move on the proof of Theorem \ref{the1.1}.

{\bf Proof of Theorem \ref{the1.1}}\quad  Let $\hat{\phi}\in \mathcal{C}^\infty_0(\mathbb{R})$ be an even, real-valued and non-negative funtion on $\mathbb{R}$ and satisfy
\begin{numcases}{\hat{\phi}(\xi)=}
1, &if $|\xi|\leq \frac{1}{4}$,\nonumber\\
0, &if $|\xi|\geq \frac{1}{2}$.\nonumber
\end{numcases}
Define the high frequency function $f_n$ and the low frequency function $g_n$ by
$$f_n=2^{-ns}\phi(x)\sin \big(\frac{17}{12}2^nx\big), \qquad g_n=\frac{12}{17}2^{-\frac{n}{k}}\phi(x), \quad n\geq1.$$
It has been showed in \cite{Li 2020} that $\|f_n\|_{B_{p, r}^\sigma}\lesssim 2^{n(\sigma-s)}$.

Set $u^n_0=f_n+g_n$, consider the system (\ref{eq5}) with initial data $u_0^n$ and $f_n$, respectively. Obviously, we have
\begin{align*}
\|u^n_0-f_n\|_{B^s_{p,r}}=\|g_n\|_{B^s_{p,r}}\leq C2^{-\frac{n}{k}},
\end{align*}
which means that
\begin{align*}
\lim_{n\to\infty}\|u^n_0-f_n\|_{B^s_{p,r}}=0.
\end{align*}
It is easy to show that
\begin{align*}
 \|u^n_0, \; f_n\|_{B^{s-1}_{p,r}}&\lesssim 2^{-\frac {n}{k}},\\
 \|u^n_0,\; f_n\|_{C^{0, 1}}&\lesssim 2^{-\frac {n}{k}}
\end{align*}
and
\begin{align*}
\|u^n_0, \; f_n\|_{B^{s+\sigma}_{p,r}}\leq C2^{\sigma n} \qquad  \mathrm{for} \qquad \sigma\geq 0,
\end{align*}
which yield
\begin{align*}
&\|u^n_0\|^{k+1}_{B^s_{p,r}}+\|u^n_0\|^k_{B^{s-1}_{p,r}}\|u^n_0\|_{B^{s+1}_{p,r}}
+\|u^n_0\|^{2k}_{B^{s-1}_{p,r}}\|u^n_0\|_{B^{s+2}_{p,r}}\lesssim 1,\\
&\|f_n\|^{k+1}_{B^s_{p,r}}+\|f_n\|^k_{B^{s-1}_{p,r}}\|f_n\|_{B^{s+1}_{p,r}}
+\|f_n\|^{2k}_{B^{s-1}_{p,r}}\|f_n\|_{B^{s+2}_{p,r}}\lesssim  1,
\end{align*}
and
\begin{align*}
&\|u^n_0\|^{k+1}_{B^s_{p,r}}+\|u^n_0\|^k_{C^{0, 1}}\|u^n_0\|_{B^{s+1}_{p,r}}
+\|u^n_0\|^{2k}_{C^{0, 1}}\|u^n_0\|_{B^{s+2}_{p,r}}\lesssim  1,\\
&\|f_n\|^{k+1}_{B^s_{p,r}}+\|f_n\|^k_{C^{0, 1}}\|f_n\|_{B^{s+1}_{p,r}}
+\|f_n\|^{2k}_{C^{0, 1}}\|f_n\|_{B^{s+2}_{p,r}}\lesssim  1.
\end{align*}
Notice that
\begin{eqnarray*}
\begin{split}
S_t(u_0^n)&=\underbrace{S_t(u_0^n)-u_0^n-t\mathbf{v}_0(u_0^n)}_{I(u_0^n)}+u_0^n+t(-(u_0^n)^k\partial_xu_0^n
+\mathbf{\mathcal{P}}(D)(u_0^n)+\mathbf{\mathcal{J}}(D)(u_0^n)),\\
S_t(f_n)&=\underbrace{S_t(f_n)-u_0^n-t\mathbf{v}_0(f_n)}_{I(f_n)}+f_n+t(-(f_n)^k\partial_xf_n
+\mathbf{\mathcal{P}}(D)(f_n)+\mathbf{\mathcal{J}}(D)(f_n)),
\end{split}
\end{eqnarray*}
then according to Proposition \ref{pro 3.3} and Lemmas \ref{lem3.2}-\ref{lem3.3}, we deduce that
\begin{align}\label{eq3.3}
\|S_{t}(u^n_0)-S_{t}(f_n)\|_{B^s_{p,r}}
\geq &~t\|(u_0^n)^k\partial_xu_0^n-(f_n)^k\partial_xf_n\|_{B_{p, r}^s}-t\|\mathbf{\mathcal{P}}(D)(u_0^n)
-\mathbf{\mathcal{P}}(D)(f_n)\|_{B_{p, r}^s}\nonumber\\
&~-t\|\mathbf{\mathcal{J}}(D)(u_0^n)
-\mathbf{\mathcal{J}}(D)(f_n)\|_{B^s_{p,r}}-\|g_n\|_{B_{p, r}^s}-Ct^2\nonumber\\
\geq&~ t\|(u_0^n)^k\partial_xu_0^n-(f_n)^k\partial_xf_n\|_{B_{p, r}^s}-C2^{-\frac {n}{k}}-Ct^{2}.
\end{align}

Moreover,
$$
(u^n_{0})^k\partial_xu^n_{0}-(f_n)^k\partial_xf_n=(g_n)^k\partial_xf_n
+(u^n_{0})^k\partial_xg_n+((u^n_{0})^k-(f_n)^k-(g_n)^k)\partial_xf_n.
$$
With the aid of Lemma \ref{lem2.1}, we find that
\begin{align*}
\|((u^n_{0})^k-f_n^k-g^k_n)\partial_xf_n\|_{B^s_{p,r}}&\leq \|(u^n_{0})^k-(f_n)^k-(g_n)^k\|_{L^\infty}\|f_n\|_{B^{s+1}_{p,r}}\\
&\;\;+\|(u^n_{0})^k-(f_n)^k-(g_n)^k\|_{B^{s}_{p,r}}\|\partial_xf_n\|_{L^\infty}\\
&\lesssim 2^{-n(s-1)},\\
\|(u^n_{0})^k\partial_xg_n\|_{B^s_{p,r}}&\leq \|u^n_0\|^k_{B^s_{p,r}}\|g_n\|_{B^{s+1}_{p,r}}\lesssim 2^{-\frac nk}.
\end{align*}
However, using the fact that
\begin{align*}
  \Delta_j((g_n)^k\partial_xf_n)=0, \;\;\; j\neq n
\end{align*}
and
\begin{align*}
  \Delta_j((g_n)^k\partial_xf_n)=(g_n)^k\partial_xf_n, \;\;\; n\geq 5,
\end{align*}
direct calculation shows that for $n\gg 1$,
\begin{align*}
  &\|(g_n)^k\partial_xf_n\|_{B_{p, r}^s}=2^{ns}\|(g_n)^k\partial_xf_n\|_{L^p}\\
  =&\|(\frac{12}{17})^k2^{-n}\phi^k\partial_x\phi\sin(\frac{17}{12}2^nx)
  +(\frac{12}{17})^{k-1}\phi^{k+1}\cos(\frac{17}{12}2^nx)\|_{L^p}\\
  \gtrsim &\|(\frac{12}{17})^{k-1}\phi^{k+1}\cos(\frac{17}{12}2^nx)\|_{L^p}-2^{-n}\rightarrow (\frac{12}{17})^{k-1} \big(\frac{\int_0^{2\pi}|\cos x|^pdx}{2\pi}\big)^{\frac 1p}\|\phi^{k+1}\|_{L^p},
\end{align*}
as $n\rightarrow\infty$ by the Riemann Theorem.

 Taking the above estimates into (\ref{eq3.3}) yields
  \begin{align*}
   \liminf_{n\rightarrow \infty}\|S_{t}(u^n_0)-S_{t}(f_n)\|_{B^s_{p,r}}\gtrsim t\quad\text{for} \ t \ \text{small enough}.
  \end{align*}
This completes the proof of Theorem \ref{the1.1}.

\section*{Acknowledgments}
 Y. Yu is supported by the Natural Science Foundation of Anhui Province (Grant No.1908085QA05) and Y. Xiao is supported by the National Natural Science Foundation of China (Grant No.11901167).

\end{document}